\documentclass[10pt,draft
]{article}

\usepackage[latin1]{inputenc}
\usepackage{amsfonts}
\usepackage{amsmath}
\usepackage{amssymb}
\usepackage{amscd}
\usepackage{amsthm}
\usepackage{indentfirst}
\usepackage[hmargin=3.6cm
,vmargin=4.0cm
]{geometry}

\usepackage{epsfig}

\newtheorem{thm}{Theorem}[section]
\newtheorem{lemma}[thm]{Lemma}
\newtheorem{prop}[thm]{Proposition}
\newtheorem{coroll}[thm]{Corollary}

\theoremstyle{definition}

\newtheorem{rem}[thm]{Remark}
\newtheorem{exam}[thm]{Example}

\newtheorem*{acknow}{Acknowledgements}
\newtheorem*{organization}{Organization of the paper}
\newtheorem*{prf}{Proof}

\newcommand{\R}{{\mathbb{R}}}

\newcommand{\Z}{{\mathbb{Z}}}

\newcommand{\cG}{{\mathcal{G}}}

\newcommand{\fc}{{:\ }}

\newcommand{\ol}{\overline}

\DeclareMathOperator{\im}{im}
\DeclareMathOperator{\id}{id}

\DeclareMathOperator{\pr}{pr}

\DeclareMathOperator{\pt}{pt}

\DeclareMathOperator{\supp}{supp}

\DeclareMathOperator{\Cont}{Cont}

\DeclareMathOperator{\osc}{osc}

\DeclareMathOperator{\Spec}{Spec}

\title{Geometric structures on contactomorphism groups and contact rigidity in jet spaces}
\author{Frol Zapolsky\footnote{Mathematisches Institut, Ludwig-Maximilian-Universit\"at, Munich, Germany; \texttt{zapolsky@math.lmu.de}}}
\begin{document}

\renewcommand{\labelenumi}{(\roman{enumi})}

\maketitle

\begin{abstract}
For a closed connected manifold $N$, we establish the existence of geometric structures on various subgroups of the contactomorphism group of the standard contact jet space $J^1N$, as well as on the group of contactomorphisms of the standard contact $T^*N \times S^1$ generated by compactly supported contact vector fields. The geometric structures are biinvariant partial orders (for $J^1N$ and $T^*N \times S^1$) and biinvariant integer-valued metrics ($T^*N\times S^1$ only). Also we prove some forms of contact rigidity in $T^*N \times S^1$, namely that certain (possibly singular) subsets of the form $X \times S^1$ cannot be disjoined from the zero section by a contact isotopy, and in addition that there are restrictions on the kind of contactomorphisms of $T^*N\times S^1$ which are products of pairwise commuting contactomorphisms generated by vector fields supported in sets of the form $U \times S^1$ with $U \subset T^*N$ Hamiltonian displaceable. The method is that of generating functions for Legendrians in jet spaces.
\end{abstract}

\section{Introduction and results}\label{section_intro}

%\subsection{Introduction and overview of the results}

In \cite{Eliashberg_Polterovich_Partially_ordered_grps_geom_cont_transfs} Eliashberg and Polterovich investigated what possible geometric structures can be put on the (universal cover of the) group of contactomorphisms of a contact manifold. They came up with the notion of a positive contact isotopy and the relation it induces on the universal cover of the identity component of the contactomorphism group. It turns out that for certain contact manifolds\footnote{Such manifolds are nowadays called orderable.} this relation is a genuine partial order; examples include standard contact projective spaces (\cite{Eliashberg_Polterovich_Partially_ordered_grps_geom_cont_transfs} using Givental's nonlinear Maslov index \cite{Givental_Nonlinear_gen_Maslov_index}), standard contact $\R^{2n+1}$ (\cite{Bhupal_partial_order_Rn}), jet spaces of closed manifolds (implicit in \cite{Colin_Ferrand_Pushkar_positive_isotopies_Leg_appl} and \cite{Chernov_Nemirovski_Leg_links_causality_Low_conj}) and their cosphere bundles (\cite{Eliashberg_Polterovich_Partially_ordered_grps_geom_cont_transfs}, \cite{Eliashberg_Kim_Polterovich_Geom_cont_transfs_dom_orderability}, \cite{Chernov_Nemirovski_nonnegative_Leg_isotopy}), and $\R^{2n} \times S^1$ (\cite{Sandon_CH_capacity_nonsqueezing_gf}).

We continue this theme in the present paper and show that there are natural biinvariant partial orders on various subgroups of the contactomorphism groups of the standard contact $J^1N = T^*N \times \R$ and $V=T^*N \times S^1$ where $N$ is a closed connected manifold of positive dimension. In order to formulate this first result precisely, let us introduce the groups with which we will be dealing. The contact structures on $J^1N$ and $V$ are given by $\xi=\ker \sigma$ where $\sigma = dz + \lambda$ is the standard contact form. Here $z$ is the coordinate on $\R$ or $S^1 = \R/\Z$ and $\lambda = p\,dq$ is the Liouville form on $T^*N$. Since the contact structure $\xi$ is fixed once and for all, we omit it from the notation throughout. We let $\Cont_0(J^1N)$ be the set of time-$1$ maps of contact isotopies of $J^1N$ supported in sets of the form $K \times \R$ where $K \subset T^*N$ is a compact subset (which may depend on the isotopy). The subgroup $\Cont_{0,c}(J^1N)$ consists of those elements of $\Cont_0(J^1N)$ whose generating isotopies are compactly supported. The subgroup $\Cont_0^\Z(J^1N) \subset \Cont_0(J^1N)$ is comprised of time-$1$ maps of contact isotopies which are equivariant with respect to the action of $\Z$ on the $\R$ factor of $J^1N = T^*N \times \R$. Finally, $\Cont_{0,c}(V) = \Cont_{0,c}(T^*N\times S^1)$ is the group of time-$1$ maps of compactly supported contact isotopies of $V$. Note that the natural projection $\Cont_0^\Z(J^1N) \to \Cont_{0,c}(V)$ is a group isomorphism, whose inverse is given by lifting. We have:
\begin{thm}\label{thm_partial_orders}The groups $\Cont_0(J^1N)$, $\Cont_{0,c}(J^1N)$, $\Cont_0^\Z(J^1N) \simeq \Cont_{0,c}(V)$ carry natural intrinsically defined biinvariant partial orders, denoted $\leq$, $\leq_c$, $\leq^\Z$, respectively. The identity maps $(\Cont_{0,c}(J^1N),\leq) \to(\Cont_{0,c}(J^1N),\leq_c)$, $(\Cont_0^\Z(J^1N),\leq) \to (\Cont_0^\Z(J^1N),\leq^\Z)$ are monotone; here $\leq$ denotes the partial orders on $\Cont_{0,c}(J^1N)$, $\Cont_0^\Z(J^1N)$ induced from the inclusion into $\Cont_0(J^1N)$.
\end{thm}
These partial orders are related to the notion of positivity introduced by Eliashberg and Polterovich in \cite{Eliashberg_Polterovich_Partially_ordered_grps_geom_cont_transfs}, in the following way. Let $(Y,\zeta)$ be a contact manifold with a cooriented contact structure $\zeta$, that is there is a contact form $\beta$ on $Y$ with $\zeta = \ker \beta$. A contact isotopy of $Y$ is called nonnegative if its generating contact vector field $X_t$ points in the nonnegative direction relative to $\zeta$, that is $\beta(X_t) \geq 0$ for all $t$. Given two elements $\widetilde\phi,\widetilde\psi$ of the universal cover $\widetilde\Cont_0(Y,\zeta)$ of the identity component of the contactomorphism group\footnote{If $Y$ is open, some restrictions on the supports of isotopies may apply.} of $Y$ we write $\widetilde\phi \leq_{\text{\rm EP}} \widetilde\psi$ if there is a nonnegative path of contactomorphisms starting at $\widetilde\phi$ and ending at $\widetilde\psi$. In \cite{Eliashberg_Polterovich_Partially_ordered_grps_geom_cont_transfs} it is shown that $\widetilde\phi \leq_{\text{\rm EP}} \widetilde\psi$ if and only if the elements $\widetilde\phi,\widetilde\psi$ are generated by contact Hamiltonians $H,K$, respectively, which satisfy $H(t,y)\leq K(t,y)$ for all $t\in[0,1]$, $y\in Y$. It follows from \cite{Colin_Ferrand_Pushkar_positive_isotopies_Leg_appl} and \cite{Chernov_Nemirovski_Leg_links_causality_Low_conj} that that the relation $\leq_{\text{\rm EP}}$ is a partial order on $\widetilde\Cont_0(J^1N)$, and in fact that it is a partial order already on $\Cont_0(J^1N)$. We have:
\begin{prop}\label{prop_EP_implies_spectral_partial_order}If $H\leq K$ are contact Hamiltonians on $J^1N$ whose time-$1$ maps are $\phi,\psi \in \Cont_0(J^1N)$, respectively, then $\phi \leq \psi$, therefore the covering projection $(\widetilde\Cont_0(J^1N),\leq_{\text{\rm EP}}) \to (\Cont_0(J^1N),\leq)$ and the identity map $(\Cont_0(J^1N),\leq_{\text{\rm EP}}) \to (\Cont_0(J^1N),\leq)$ are monotone.
\end{prop}
\noindent Analogous statements are valid for the other groups introduced above, together with their partial orders; see subsection \ref{section_partial_orders} below.

The next result concerns the existence of biinvariant $\Z$-valued metrics on $\Cont_{0,c}(V)$. In \cite{Sandon_integer_biinvt_metric_Rtwo_n_S_one} Sandon constructed a biinvariant $\Z$-valued metric on the identity component of the contactomorphism group of the standard contact $\R^{2n}\times S^1$. 
Denote by $\preceq$ the partial order on $\Cont_{0,c}(V)$ induced from the isomorphism $\Cont_{0,c}(V) \simeq \Cont_0^\Z(J^1N)$. We have
\begin{thm}\label{thm_metrics}
There are two natural equivalent biinvariant metrics on $\Cont_{0,c}(V)$, which we denote $\rho_{\osc}$ and $\rho_{\sup}$. The triple $(\Cont_{0,c}(V),\preceq,\rho_{\sup})$ is an ordered metric space in the sense that $\phi\preceq\chi\preceq\psi$ implies $\rho_{\sup}(\phi,\chi)\leq \rho_{\sup}(\phi,\psi)$. There is an order-preserving isometric group embedding of $(\Z,\leq,|\cdot-\cdot|)$ into this ordered metric group, in particular it has infinite diameter.
\end{thm}
\begin{rem}Intuitively speaking, the metric $\rho_{\osc}$ is analogous to the oscillation metric on functions, whereas $\rho_{\sup}$ is analogous to the uniform metric. To further amplify this analogy, note that for functions $f,g,h$ the condition $f\leq g \leq h$ implies $\|f-g\|_{C^0} \leq \|f-h\|_{C^0}$, but this conclusion fails if the uniform metric is replaced with oscillation.
\end{rem}

\begin{rem}There is work in progress by Fraser and Polterovich \cite{Fraser_Polterovich_Sandon_type_metrics_cont_groups} aimed, among other things, at showing the existence of $\Z$-valued biinvariant metrics on contactomorphism groups of a class of circle bundles, which includes $T^*N \times S^1$.
\end{rem}

This concludes the first part of the results. The second part concerns contact rigidity in $T^*N \times S^1$. To put it into context, let us mention the well-known phenomenon of symplectic rigidity, which manifests itself, for example, in the fact that certain subsets of symplectic manifolds cannot be disjoined from other subsets by a Hamiltonian isotopy, while there may be no topological restrictions for this. An example is the zero section $O$ of $T^*N$, for which symplectic rigidity is a version of the Arnold conjecture. By now it is known that there are also singular subsets which exhibit symplectic rigidity, see \cite{Entov_Polterovich_rigid_subsets_sympl_mfds}. For example, recall the partial symplectic quasi-state $\eta \fc C^\infty_c(T^*N) \to \R$ constructed in \cite{Monzner_Vichery_Zapolsky_partial_qms_qss_cot_bundles} (there it is denoted $\zeta$ or $\zeta_0$). Following \cite{Entov_Polterovich_rigid_subsets_sympl_mfds}, call a compact subset $X \subset T^*N$ superheavy if whenever $f \in C^\infty_c(T^*N)$ is such that $f|_X = c \in \R$, then $\eta(f) = c$. One of the consequences of this definition is the following result:
\begin{thm}[\cite{Monzner_Vichery_Zapolsky_partial_qms_qss_cot_bundles}]A superheavy subset cannot be displaced from the zero section by a Hamiltonian isotopy.
\end{thm}
\noindent For the convenience of the reader, let us give examples of superheavy subsets. Recall that a subset $Z$ of a symplectic manifold is called displaceable if there is a Hamiltonian diffeomorphism $\phi$ of the symplectic manifold with $\phi(Z) \cap \ol Z = \varnothing$.
\begin{exam}[\cite{Monzner_Vichery_Zapolsky_partial_qms_qss_cot_bundles}]
\begin{enumerate}
\item The zero section is superheavy;
\item $X$ is superheavy if $T^*N - X = U_\infty \cup \bigcup_j U_j$ is a finite disjoint union, where $U_\infty$ is the unbounded connected component of the complement of $X$ (or the union of the two components when $N = S^1$), such that $U_\infty \cap O = \varnothing$, while $U_j$ all have the property that $\eta|_{C^\infty_c(U_j)}\equiv 0$, which happens, for instance, when every compact contained in $U_j$ is displaceable;
\item if $X_1,\dots,X_k$ are superheavy in $T^*N_1,\dots,T^*N_k$, respectively, then the product $\prod_jX_j$ is superheavy in $T^*\prod_jN_j$;
\item if $X$ is superheavy, so is $\phi(X)$ for a Hamiltonian diffeomorphism $\phi(X)$.
\end{enumerate}
\end{exam}
\noindent There is an analogous phenomenon in contact topology, which takes the form of another version of the Arnold conjecture: the zero section $O \subset J^1N$ (or $O\subset V$) cannot be disjoined from the zero wall $O \times \R \subset J^1N$ (or $O \times S^1 \subset V$) by a contact isotopy. Since the zero section in $T^*N$ is superheavy, our next result subsumes both these phenomena, namely we have
\begin{thm}\label{thm_rigidity_superheavy} Let $X \subset T^*N$ be a superheavy subset. Then for any $\alpha \in \Cont_{0,c}(V)$ we have $\alpha(X\times S^1) \cap O \neq \varnothing$.
\end{thm}

Finally, we establish another form of contact rigidity, namely we prove:
\begin{thm}\label{thm_rigidity_products} Let $U_1,\dots,U_k \subset T^*N$ be displaceable open subsets and $\phi_j \in \Cont_{0,c}(U_j \times S^1)$ for $j=1,\dots,k$. Assume in addition that the $\phi_j$ all pairwise commute. Let $H$ be a compactly supported contact Hamiltonian on $V$ whose time-$1$ map is the product $\phi_1\dots\phi_k$. Then $\min_{t\in[0,1],y\in O\times S^1}|H(t,y)| = 0$. In other words, the time-$1$ map of a contact Hamiltonian which is bounded away from zero along the zero wall cannot be the product of commuting contactomorphisms as above.
\end{thm}
\noindent A similar phenomenon is observed in the symplectic case, see \cite{Entov_Polterovich_Zapolsky_qms_Poisson_bracket}, \cite{Monzner_Vichery_Zapolsky_partial_qms_qss_cot_bundles}.

\begin{organization}The next subsection comments on the method used to obtain the results, while in subsection \ref{section_prelims_notations} we give some preliminaries on contact and symplectic geometry, as well as fix the notation. In section \ref{section_sp_numbers_contact} we define spectral numbers for Legendrians in $J^1N$ and for contactomorphisms of $J^1N$ and $T^*N\times S^1$ and prove their properties. Finally, section \ref{section_proofs} contains the proofs of the main results.
\end{organization}

\begin{acknow}I wish to thank Sheila ``Margherita'' Sandon and Leonid Polterovich for valuable discussions, and Bijan Sahamie for listening to a preliminary version of the above results. This work was carried out at Ludwig-Maximilian-Universit\"at, Munich, and I would like to acknowledge its excellent research atmosphere and hospitality.
\end{acknow}

\subsection{Method}\label{section_method}

All of the above results are consequences of the existence and uniqueness of generating functions quadratic at infinity for Legendrian submanifolds of $J^1N$ which are isotopic to the zero section through Legendrian submanifolds. Details are given in section \ref{section_sp_numbers_contact}. For now let us mention that there are certain numbers, called spectral and denoted $\ell(A,L)$, attached to Legendrians $L \subset J^1N$ and homology classes $A \in H_*(N)$. We let $\ell_+ = \ell([N],\cdot)$ and $\ell_- = \ell(\pt,\cdot)$. Next, we define $\ell_\pm \fc \Cont_0(J^1N) \to \R$ as $\ell_\pm(\phi)=\ell_\pm\big(\phi(O)\big)$.

The partial order $\leq$ is given by declaring that $\phi\leq\psi$ if $\ell_+(\alpha\phi\psi^{-1}\alpha^{-1})\leq 0$ for all $\alpha \in \Cont_0(J^1N)$. Similar definitions hold for the other partial orders appearing in theorem \ref{thm_partial_orders}.

To define the metrics on $\Cont_{0,c}(V)$, let $\phi \mapsto \widetilde\phi$ denote the inverse of the isomorphism $\Cont_0^\Z(J^1N) \to \Cont_{0,c}(V)$ given by projection. Define $\ell_\pm \fc \Cont_{0,c}(V) \to \R$ via $\ell_\pm(\phi) = \ell_\pm\big(\widetilde\phi\big)$. The metrics $\rho_{\osc}$ and $\rho_{\sup}$ are given by
$$\rho_{\osc}(\phi,\psi) = \max\big\{\lceil \ell_+(\alpha\phi\psi^{-1}\alpha^{-1})\rceil - \lfloor\ell_-(\alpha\phi\psi^{-1}\alpha^{-1})\rfloor\,|\,\alpha\in\Cont_{0,c}(V)\big\}$$
and
$$\rho_{\sup}(\phi,\psi) = \max\big\{\max\big(|\lceil \ell_+(\alpha\phi\psi^{-1}\alpha^{-1})\rceil|,|\lfloor\ell_-(\alpha\phi\psi^{-1}\alpha^{-1})\rfloor|\big)\,|\,\alpha\in\Cont_{0,c}(V)\big\}\,.$$

Finally, contact rigidity is a consequence of the existence of a certain ``partial quasi-morphism'' $\nu \fc \Cont_{0,c}(V) \to \R$. It is proved below that $\lceil\ell_+(\phi\psi)\rceil \leq \lceil\ell_+(\phi)\rceil + \lceil\ell_+(\psi)\rceil$, which means that the sequence $\lceil\ell_+(\phi^k)\rceil$ is subadditive (with respect to $k$). Therefore we can define
$$\nu(\phi)=\lim_{k\to\infty}\frac{\ell_+(\phi^k)}{k}\,.$$
In subsection \ref{section_cont_rigidity} we prove the properties of this function which imply the above rigidity results.

\subsection{Preliminaries and notations}\label{section_prelims_notations}

The cotangent bundle $T^*N$ admits a canonical symplectic form $\omega = d\lambda = dp \wedge dq$. A Hamiltonian $h$ on $T^*N$ is a compactly supported smooth function on $[0,1]\times T^*N$; its time-dependent Hamiltonian vector field is defined via $\omega(X_{h_t},\cdot)=-dh_t$, and its Hamiltonian flow $\phi_h^t$ is obtained by integrating $X_{h_t}$. We let $\phi_h = \phi_h^1$.

A contact Hamiltonian (or just a Hamiltonian if no confusion can arise) on $J^1N$ or $V$ is a smooth time-dependent function, which in the case of $V$ is usually assumed to have compact support. The time-dependent contact vector field $X_{H_t}$ of a contact Hamiltonian $H_t$ is defined by requiring $\alpha(X_{H_t}) = H_t$ and $d\alpha(X_{H_t},\cdot)=dH_t(R_\alpha)\alpha - dH_t$, where $R_\alpha = \partial_z$ is the Reeb vector field of $\alpha$. The contact flow of $H$ is the flow $\phi_H^t$ obtained by integrating $X_{H_t}$. We put $\phi_H = \phi_H^1$.

The connection between the two geometries is as follows. Let $\pi \fc J^1N \to T^*N$ or $\pi \fc V=T^*N \times S^1 \to T^*N$ be the symplectic projection which forgets the $\R$ or $S^1$ factor. If $h$ is a Hamiltonian on $T^*N$ then $H=\pi^*h$ is a contact Hamiltonian on $J^1N$ or $V$. The contact vector field of $H$ satisfies $\pi_*X_H = X_h$, in particular, the projection $\pi$ intertwines the contact and the Hamiltonian flows, that is $\pi\circ\phi_H^t = \phi_h^t\circ\pi$.

Finally, for a function $f \in C^\infty(M)$ we let $j^1f = \{(q,-d_qf,f(q))\in J^1N\,|\,q\in N\} \subset J^1N$ be its $1$-jet. It is a Legendrian submanifold.

\section{Spectral numbers for contactomorphisms}\label{section_sp_numbers_contact}

In this section we define the spectral numbers on the above contactomorphism groups and prove the relevant properties which will be used in the proofs of the main results. All homology is with coefficients in $\Z_2$.

\subsection{Spectral numbers for Legendrians in jet spaces via generating functions}\label{section_sp_numbers_Leg}

Here we describe the spectral numbers for Legendrian submanifolds of the jet space $J^1N$. In order to define these, we need the notion of a generating function quadratic at infinity. We start with the more general notion of a generating function for a (possibly singular) Legendrian submanifold of $J^1N$. For more information see \cite{Viterbo_gfqi}, \cite{Chaperon_On_generating_families}, \cite{Theret_PhD_thesis}, \cite{Chekanov_critical_pts_quasifcns_gf_Leg}.

If $S\fc N \times E \to \R$ is a smooth function, where $E$ is a finite-dimensional vector space, we put $\Sigma_S = \{(q,e) \in N \times E\,|\,\partial_eS(q,e)=0\}$ and $i_S \fc \Sigma_S \to J^1N$ is defined by $i_S(q,e) = (q,-\partial_qS(q,e),S(q,e))$. The subset $\Sigma_S$ carries a canonical tangent distribution, namely the kernel of the differential of the fiberwise derivative map $N \times E \to E^*$, $(q,e) \mapsto \partial_eS$, at points mapped to $0 \in E^*$, and in this sense the map $i_S$ has a differential defined on this tangent distribution, which satisfies $(i_S)^*\sigma = 0$, which means that its image $i_S(\Sigma_S)$ is in a canonical sense Legendrian (even if it fails to be a submanifold). Note that the critical points of $S$ are mapped by $i_S$ to the zero wall $O\times \R$ and the height of the image in $\R$ of a critical points equals its critical value. We call the set of critical values of $S$ its spectrum, denoted $\Spec(S) \subset \R$. We also call this subset the spectrum of the Legendrian $L=i_S(\Sigma_S)$ and denote it by $\Spec(L)$.

If the fiberwise derivative is transverse to $0 \in E^*$ then $\Sigma_S$, being the preimage of a regular value, is a submanifold whose tangent distribution is the one described above, and its image $L=i_S(\Sigma_S)$ is an immersed Legendrian submanifold of $J^1N$. We say that $S$ generates $L$. The projection $\pi(L) \subset T^*N$ is an immersed Lagrangian submanifold, and we also say that $S$ generates it.

There are certain operations on generating functions which will be used below. Namely, if $S \fc N \times E \to \R$ and $S' \fc N' \times E' \to \R$ are generating functions, their fiberwise direct sum and difference $S \oplus S',\,S\ominus S' \fc N \times E\times E' \to \R$ are defined by $(S \oplus S')(q,e,e')= S(q,e) + S'(q,e')$ and similarly for the difference. Given two subsets $L,L' \subset J^1N$ their sum and difference are given by
$$L\pm L' = \{(q,p\pm p',z \pm z')\,|\,(q,p,z)\in L,\,(q,p',e')\in L'\}\,;$$
we let $-L:=O-L$ where $O$ is the zero section. It follows that if $S$ generates $L$ and $S'$ generates $L'$, then $S\oplus S'$, $S \ominus S'$, and $-S$ generate $L + L'$, $L-L'$, and $-L$, respectively. Note that even when $L,L'$ are submanifolds, $L\pm L'$ need not be.

\begin{rem}Usually the Legendrian generated by $S$ is given by $\ol L = \ol i_S(\Sigma_S)$, where $\bar i_S\fc\Sigma_S \to J^1N$, $(q,e) \mapsto (q,\partial_qS,S(q))$. This has to do with the fact that this formula assumes that the contact structure on $J^1N$ is given by $\ker(dz - \lambda)$, while we use $\ker(dz + \lambda)$. If we let $\iota \fc (J^1N,\ker(dz - \lambda)) \to (J^1N,\ker(dz+\lambda))$ be the contactomorphism $\iota(q,p,e)=(q,-p,e)$, then the Legendrians $\ol L \subset (J^1N,\ker(dz-\lambda))$ and $L \subset (J^1N,\ker(dz+\lambda))$ are connected  via $L = \iota(\ol L)$.
\end{rem}

We say that $S$ is quadratic at infinity, abbreviated gfqi, if there is a nondegenerate quadratic form $Q \fc E \to \R$ such that for points $(q,e)$ outside a compact subset of $N\times E$ we have $S(q,e) = Q(e)$. The following is the cornerstone of the whole story:
\begin{thm}[\cite{Chaperon_On_generating_families}, \cite{Theret_PhD_thesis}, \cite{Chekanov_critical_pts_quasifcns_gf_Leg}]Let $L \subset J^1N$ be a Legendrian submanifold which is isotopic through Legendrian submanifolds to the zero section. Then it is generated by a gfqi. Moreover if $L_t$ is an isotopy consisting of Legendrian submanifolds Legendrian isotopic to the zero section, then there is a smooth family $S_t \fc N \times E \to \R$ of gfqi such that $S_t$ generates $L_t$ for each $t$. These generating functions are unique up to gauge transformation and stabilization.
\end{thm}
\begin{rem}A gauge transformation is a compactly supported diffeomorphism of $N \times E$ which preserves the fibers $\{q\}\times E$. A stabilization of a gfqi $S$ is the gfqi $S\oplus Q' \fc N \times E \times E' \to \R$, where $E'$ is a finite-dimensional vector space, $Q'$ a nondegenerate quadratic form on it and $(S\oplus Q')(q,e,e')=S(q,e)+Q'(e')$. These two operations, applied to a gfqi, keep intact the Legendrian it generates. They also preserve the Morse-theoretic invariants obtained from a gfqi \cite{Viterbo_gfqi}.
\end{rem}

This allows us to define spectral numbers for Legendrians in a standard way. Namely, let $S \fc N \times E \to \R$ be a gfqi; then for any $a,b \in \R$ such that $a>b$ and $b$ is sufficiently negative, the homotopy type of the pair $(\{S<a\},\{S<b\})$ does not depend on $b$. Moreover for $a$ large enough the homotopy type of this pair is also independent of $a$ and in fact its homology is canonically graded isomorphic to $H_*(N\times E_-,N\times (E_- - \{0\}))$, where $E_- \subset E$ is the negative subspace of the quadratic form associated to $S$. This latter homology, with coefficients in $\Z_2$, is isomorphic, via the K\"unneth isomorphism, to $H_*(N) \otimes \Z_2$, where the $\Z_2$ lives in degree $\dim E_-$. For $b$ small enough and a homology class $A \in H_*(N;\Z_2)$ let
$$\ell(A,S)=\inf\{a\in\R\,|\,A\otimes 1 \in \im i^a\}\,,$$
where $i^a \fc H_*(\{S<a\},\{S<b\}) \to H_*(N) \times \Z_2$ is the composition of the map induced on homology by inclusion and the K\"unneth isomorphism, and $1 \in \Z_2$ is the generator. Note that by Lusternik-Schnirelman theory $\ell(A,S)$ is a critical value of $S$, thus $\ell(A,S) \in \Spec(S)$. If $S$ generates a Legendrian $L$, we set $\ell(A,L):=\ell(A,S)$. Since gauge transformations and stabilizations do not alter the spectral invariants of a gfqi, it follows that $\ell(A,L)$ is well-defined, that is it only depends on $L$.

We now indicate the connection of these Legendrian spectral numbers to Lagrangian spectral invariants defined on the Hamiltonian group of $T^*N$, as described in \cite{Monzner_Vichery_Zapolsky_partial_qms_qss_cot_bundles}, for example. If $\phi$ is a Hamiltonian diffeomorphism of $T^*N$, let $\ell^{T^*N}(A,\phi)$, where $A \in H_*(N)$, be the spectral invariants appearing \emph{ibid.} Since $\phi(O)$ is a Lagrangian which is exactly isotopic to the zero section, it admits a gfqi $S$, which moreover can be normalized so that its spectrum $\Spec(S)$ coincides with the action spectrum $\Spec(\phi)$ (see \cite{Monzner_Vichery_Zapolsky_partial_qms_qss_cot_bundles} for a definition). When we normalize $S$ in this way, we obtain (see \cite{Monzner_Vichery_Zapolsky_partial_qms_qss_cot_bundles}):
$$\ell^{T^*N}(A,\phi) = \ell(A,S)\,.$$
This gfqi $S$ also generates a Legendrian lift $L$ of $\phi(O)$. This Legendrian may be obtained from $\phi$ as follows. Let $h$ be a Hamiltonian on $T^*N$ such that $\phi=\phi_h$ and let $H = \pi^*h$ be the corresponding contact Hamiltonian on $J^1N$. Then $L=\phi_H(O)$. In conclusion of this discussion, we have
$$\ell^{T^*N}(A,\phi_h) = \ell(A,\phi_{\pi^*h}(O))\,.$$

We next give the properties of these spectral numbers which will be needed in the sequel.
\begin{prop}\label{prop_properties_sp_nums_Leg}
Let $L \subset J^1N$ be an embedded Legendrian submanifold which is isotopic to the zero section through embedded Legendrian submanifolds. Then to any homology class $A \in H_*(N)$ there is associated a number $\ell(A,L)$ with the following properties:
\begin{enumerate}
\item $\ell(A,L) \in \Spec(L)$;
\item $\ell(A \cap B,L+L') \leq \ell(A,L) + \ell(B,L')$, where $\cap$ is the intersection product on homology and $L'$ is another Legendrian isotopic to the zero section through Legendrians;
\item $\ell([N],L) = -\ell(\pt,-L)$;
\item if $\ell([N],L) = \ell(\pt,L) = a \in \R$ then $L = j^1a$, where $a$ is considered as a constant function on $N$;
\item for $f \in C^\infty(N)$ the numbers $\ell(A,j^1f)$ are the usual homological spectral invariants of $f$;
\item if $h$ is a Hamiltonian on $T^*N$, then $\ell^{T^*N}(A,\phi_h) = \ell(A,\phi_{\pi^*h}(O))$;
\item if $\phi_t$ is a contact isotopy of $J^1N$ then the numbers $\ell(a,\phi(L))$ and $\ell(A,L - \phi^{-1}(O))$ are simultaneously positive, negative, or zero; if $\phi_t$ is a contact isotopy of $J^1N$ such that all $\phi_t$ are equivariant with respect to the natural $\Z$ action, then either both $\ell(A,\phi(L))$ and $\ell(A, L - \phi^{-1}(O))$ are integer, in which case they are equal, or they belong to the same interval $(k,k+1)$ for some $k \in \Z$; in particular we always have
$$\lceil\ell(A,\phi(L))\rceil = \lceil\ell(A,L -\phi^{-1}(O))\rceil \quad\text{and}\quad \lfloor\ell(A,\phi(L))\rfloor = \lfloor\ell(A,L -\phi^{-1}(O))\rfloor\,;$$
\item $\ell(A,L-L) = 0$.
\end{enumerate}
\end{prop}

\begin{prf}
(i) follows directly from definition.

(ii-iii) are proved in \cite{Viterbo_gfqi}.

(iv) If $S$ is a gfqi, it restricts to a gfqi $S_q = S|_{\{q\}\times E} \fc E \to \R$. This function has only one spectral invariant, call it $\ell_q(S)$. The comparison inequality (\cite{Viterbo_gfqi}) says $\ell_-(S) \leq \ell_q(S) \leq \ell_+(S)$ for all $q$, which means that $\ell_q(S) = a$ for all $q$. Now in general the function $q \mapsto \ell_q(S)$ is Lipschitz and also smooth on an open dense subset of $N$. Its $1$-jet, defined on this open subset, has as its image a part of $L$ which is graphical over $N$. In our case this means that $L$ contains $j^1a$, since $\ell_q(S) = a$, therefore $L = j^1a$.

(v) In this case $f$ is a gfqi for its jet and the claim follows.

(vi) follows from the discussion above.

(vii) The first assertion is proved by Bhupal \cite{Bhupal_partial_order_Rn}, the second one by Sandon \cite{Sandon_CH_capacity_nonsqueezing_gf}. In fact, there she proves that if the function $t\mapsto \ell(A,\phi_t^{-1}\phi(L) - \phi_t^{-1}(O))$ attains an integer value then it is constant, and of course equal to that value. The statement follows upon taking $t=0,1$.

(viii) follows from the previous point. Indeed, if $L = \phi(O)$ where $\phi$ is the time-$1$ map of a contact isotopy then (vii) says that $\ell(A,L-L) = \ell(A,\phi(O) - \phi(O)) =0 $ if and only if $\ell(A,\phi^{-1}\phi(O)) = 0$ but the latter number is $\ell(A,O)=0$. \qed
\end{prf}

We will also need to use the dependence of spectral numbers of a Legendrian $L$ on the contact Hamiltonian generating an isotopy which maps $O$ to $L$. Let $H$ be a contact Hamiltonian generating the contact isotopy $\phi_t$ on $J^1N$; fix two Legendrians $L,L_0$ and let $L_t = \phi_t(L_0)$. We will need the following fact. The reader will find more details in \cite{Colin_Ferrand_Pushkar_positive_isotopies_Leg_appl}.
\begin{prop}
Let $L$ have a gfqi $S$ and $L_t$ a gfqi $S_t$. Consider the Cerf diagram of the gfqi $S_t \ominus S$, that is the subset $\{(t,c)\,|\,t\in[0,1],\,c\in\Spec(S_t \ominus S)\} \subset [0,1]\times \R$. Then its slope at a smooth point $(t,c)$ represented by a point $(q,p,e) \in L_t$ equals $H(q,p,e)$. \qed
\end{prop}
\noindent Since $S_t\ominus S$ is a gfqi generating $L_t - L$, the following follows from the continuity (see \cite{Viterbo_gfqi}) and spectrality (point (i) of proposition \ref{prop_properties_sp_nums_Leg}) of spectral numbers:
\begin{lemma}\label{lemma_sp_nums_and_ct_Hamiltonian} Let $L,L_0\subset J^1N$ be two Legendrians isotopic through Legendrians to the zero section and let $H_t$ be a time-dependent contact Hamiltonian. Put $L_t = \phi_H^t(L_0)$. Then
$$\int_0^1 \min_{L_t} H_t \,dt \leq \ell(A,L_1-L) - \ell(A,L_0-L) \leq \int_0^1 \max_{L_t} H_t\,dt\,;$$
in particular,
$$\int_0^1 \min_{L_t} H_t \,dt \leq \ell(A,L_1-L_0) \leq \int_0^1 \max_{L_t} H_t\,dt\,.$$
\end{lemma}
\begin{prf}
The first assertion follows from the above discussion. For the second one we only need to note that $\ell(A,L_0-L_0)=0$, which is point (viii) in proposition \ref{prop_properties_sp_nums_Leg}. \qed
\end{prf}

In the sequel we will need the particular case when $L=L_0$ and it is obtained from the zero section via a Legendrian isotopy. We will also need the following corollary.
\begin{coroll}\label{coroll_comparison}Let $H,K$ be contact Hamiltonians on $J^1N$ with $H \leq K$ everywhere. Then $\ell(A,\phi_H^t(L)) \leq \ell(A,\phi_K^t(L))$ for all $t$ and $L$ Legendrian isotopic to the zero section.
\end{coroll}
\begin{prf}
Since $H \leq K$ everywhere, it follows that for every $t$ there is a nonnegative contact isotopy beginning at $\phi_H^t$ and ending at $\phi_K^t$ (\cite{Eliashberg_Polterovich_Partially_ordered_grps_geom_cont_transfs}). Now apply lemma \ref{lemma_sp_nums_and_ct_Hamiltonian}. \qed
\end{prf}

\subsection{Spectral numbers for contactomorphisms}\label{section_subsect_sp_nums_contact}

Here we construct the spectral numbers for contactomorphisms of $J^1N$ and $V = T^*N \times S^1$. Given an element $\phi \in \Cont_0(J^1N)$ and a homology class $A \in H_*(N)$ we let $\ell(A,\phi) = \ell(A,\phi(O))$. We will only use the spectral numbers $\ell_+ = \ell([N],\cdot)$ and $\ell_- = \ell(\pt,\cdot)$. These functions have the following properties:
\begin{prop}\label{prop_properties_sp_nums_contact}
\begin{enumerate}
\item the numbers $\ell_+(\phi)$ and $-\ell_-(\phi^{-1})$ are either both zero, or both positive, or both negative;
\item if $\ell_+(\phi) \leq 0$ and $\ell_+(\psi) \leq 0$ then $\ell_+(\phi\psi) \leq 0$;
\item let $h$ be a Hamiltonian on $T^*N$, then $\ell(A,\phi_{\pi^*h}) = \ell^{T^*N}(A,\phi_h)$.
\end{enumerate}
\end{prop}
\begin{prf}
(i) We have by point (vii) of proposition \ref{prop_properties_sp_nums_Leg}: $\ell(A,\phi(L))$ and $\ell(A,L-\phi^{-1}(O))$ are always simultaneously either zero or positive or negative. Thus the same is true of $\ell_+(\phi)$ and $\ell_+(-\phi^{-1}(O)) = -\ell_-(\phi^{-1})$, which is what was claimed.

(ii) We have: $\ell_+(\psi(O) - \phi^{-1}(O)) \leq \ell_+(\psi(O)) + \ell_+(-\phi^{-1}(O)) = \ell_+(\psi) - \ell_-(\phi^{-1})$. Since by assumption $\ell_+(\phi) \leq 0$, we see that $-\ell_-(\phi^{-1}) \leq 0$, therefore $\ell_+(\psi(O) - \phi^{-1}(O)) \leq 0$, which means that $\ell_+(\phi\psi) = \ell_+(\phi\psi(O)) \leq 0$.

(iii) follows from the definition and point (vi) in proposition \ref{prop_properties_sp_nums_Leg} above. \qed
\end{prf}

The above definition, of course, specializes to the subgroup $\Cont_0^\Z(J^1N) \subset \Cont_0(J^1N)$. Also, let $\phi \mapsto \widetilde\phi$ denote the inverse of the isomorphism $\Cont_0^\Z(J^1N) \to \Cont_{0,c}(V)$, and put $\ell(A,\phi) = \ell\big(A,\widetilde\phi\big)$ for $\phi \in \Cont_{0,c}(V)$. Then
\begin{prop}\label{prop_properties_sp_nums_contact_periodic} Let $\ell(A,\cdot)$ be the functions defined on $\Cont_0^\Z(J^1N)$ and $\Cont_{0,c}(V)$. Then
\begin{enumerate}
\item (duality) $\lceil\ell_+(\phi)\rceil = -\lfloor\ell_-(\phi^{-1})\rfloor$;
\item (triangle inequalities) we have $\lceil\ell_+(\phi\psi)\rceil \leq \lceil\ell_+(\phi)\rceil + \lceil\ell_+(\phi\psi)\rceil$, $\lfloor\ell_-(\phi\psi)\rfloor \geq \lfloor\ell_-(\phi)\rfloor + \lfloor\ell_-(\psi)\rfloor$ and $\lfloor\ell_-(\phi\psi)\rfloor \leq \lfloor\ell_-(\phi)\rfloor + \lceil\ell_+(\psi)\rceil$, $\lfloor\ell_-(\phi\psi)\rfloor \leq \lceil\ell_+(\phi)\rceil + \lfloor\ell_-(\psi)\rfloor$;
\item let $H$ be either a contact Hamiltonian on $J^1N$, supported in $\text{\rm compact}\times \R$ and $\Z$-invariant or a compactly supported contact Hamiltonian on $V$, and $\alpha$ an element in either $\Cont_0^\Z(J^1N)$ or $\Cont_{0,c}(V)$; then
$$\left\lfloor\int_0^1\min H_t\,dt \right\rfloor\leq \ell_\pm(\alpha\phi_H\alpha^{-1}) \leq \left\lceil\int_0^1\max H_t\,dt \right\rceil\;.$$
\end{enumerate}
\end{prop}
\begin{prf}
(i) Point (vii) of proposition \ref{prop_properties_sp_nums_Leg} implies $\lceil \ell_+(\phi(O)) \rceil = \lceil \ell_+(O-\phi^{-1}(O)) \rceil$, while point (iii) says $\ell_+(-\phi^{-1}(O)) = -\ell_-(\phi^{-1}(O))$, therefore
$$\lceil \ell_+(\phi(O)) \rceil = \lceil -\ell_-(\phi^{-1}(O)) \rceil = -\lfloor \ell_-(\phi^{-1}(O)) \rfloor\,.$$

(ii) We only prove the first inequality, the rest being proved similarly. We have
$$\ell_+(\psi(O) - \phi^{-1}(O)) \leq \ell_+(\psi(O)) + \ell_+(-\phi^{-1}(O)) = \ell_+(\psi) - \ell_-(\phi^{-1})\,.$$
Therefore
$$\lceil \ell_+(\phi\psi) \rceil = \lceil \ell_+(\phi\psi(O)) \rceil = \lceil \ell_+(\psi(O) - \phi^{-1}(O)) \rceil \leq \lceil \ell_+(\psi) - \ell_-(\phi^{-1}) \rceil\,.$$
Since
$$\lceil \ell_+(\psi) - \ell_-(\phi^{-1}) \rceil \leq \lceil \ell_+(\psi) \rceil + \lceil - \ell_-(\phi^{-1}) \rceil = \lceil \ell_+(\psi) \rceil + \lceil \ell_+(\phi) \rceil\,,$$
we are done.

(iii) Abbreviate $\phi=\phi_H$. We have $\ell_+(\alpha\phi\alpha^{-1}) = \ell_+(\alpha\phi\alpha^{-1}(O))$. Since $\alpha \in \Cont_0^\Z(J^1N)$ by assumption, we obtain $\lceil \ell_+(\alpha\phi\alpha^{-1}(O)) \rceil = \lceil \ell_+(\phi\alpha^{-1}(O) - \alpha^{-1}(O)) \rceil$. By lemma \ref{lemma_sp_nums_and_ct_Hamiltonian} it is true that
$$\int_0^1 \min H_t \,dt \leq \ell_+(\phi\alpha^{-1}(O) - \alpha^{-1}(O))\leq \int_0^1 \max H_t \,dt\,.$$
Arguing similarly with $\ell_-(\alpha\phi\alpha^{-1})$ and taking the upper and lower integer parts, we obtain the desired inequalities. \qed
\end{prf}

For future use we record the following
\begin{lemma}\label{lemma_contactomorphism_fixing_jets_identity}
If $\phi_t$ is a contact isotopy of $J^1N$ supported in $K \times \R$, where $K \subset T^*N$ is a compact subset, such that $\phi(j^1f)=j^1(f + a(f))$ for every $f\in C^\infty(N)$, where $a \fc C^\infty(N) \to \R$ is a function, then $\phi = \id$.
\end{lemma}
\begin{prf}
Let $y_0 \in J^1N$. Since $\phi$ equals the identity map outside its support, whenever $y \notin \supp \phi$ and $f$ is such that $y,y_0 \in j^1f$, we have $\phi(y) = y$, which implies $\phi(j^1f)=j^1f$, and therefore $\phi(y_0) \in \phi(j^1f) = j^1f$. Thus
$$\phi(y_0) \in \bigcap\{j^1f\,|\,j^1f\ni y_0\text{ and } j^1f \text{ meets the complement of } \supp \phi\}\,.$$
Since for any $y \neq y_0$ there is a function $f$ with $y \notin j^1f$ and such that $j^1f$ meets the complement of $\supp \phi$, we see that the above intersection contains only one point, namely $y_0$ which implies that $\phi(y_0) = y_0$. \qed
\end{prf}

\section{Proofs}\label{section_proofs}

\subsection{Partial orders}\label{section_partial_orders}

Here we prove theorem \ref{thm_partial_orders}.
\begin{prf}[of theorem \ref{thm_partial_orders}]Let us first define the partial orders. The proof for all the three groups is identical, so let $\cG$ denote either one of them, and let $\leq$ be the binary relation on $\cG$ defined as follows: for $\phi,\psi \in \cG$ we let $\phi \leq \psi$ if $\ell_+(\alpha\phi\psi^{-1}\alpha^{-1}) \leq 0$ for all $\alpha \in \cG$. We need to show that this relation is biinvariant, reflexive, antisymmetric, and transitive.

The fact that this relation is biinvariant, follows immediately from the definition. Reflexivity is obvious. For transitivity it is enough to show that if $\phi \leq \id$ and $\psi \leq \id$ then $\phi \psi \leq \id$. By assumption, for all $\alpha \in \cG$ we have $\ell_+(\alpha\phi\alpha^{-1}),\ell_+(\alpha\psi\alpha^{-1}) \leq 0$. It follows from proposition \ref{prop_properties_sp_nums_contact} above that $\ell_+(\alpha\phi\psi\alpha^{-1}) = \ell_+(\alpha\phi\alpha^{-1}\cdot\alpha\phi\alpha^{-1}) \leq 0$, which proves what we want. In order to prove antisymmetry it suffices to show that $\phi \leq \id$ and $\phi^{-1}\leq \id$ imply $\phi = \id$. Since $\ell_+(\alpha\phi^{-1}\alpha^{-1}) \leq 0$, it follows that $-\ell_-(\alpha\phi\alpha^{-1}) \leq 0$ (point (i) of proposition \ref{prop_properties_sp_nums_contact}), which implies the following string of inequalities:
$$0\leq\ell_-(\alpha\phi\alpha^{-1})\leq \ell_+(\alpha\phi\alpha^{-1}) \leq 0\,,$$
forcing $\ell_\pm(\alpha\phi\alpha^{-1}) = 0$ (point (iv) of proposition \ref{prop_properties_sp_nums_Leg}), which means that $\alpha\phi\alpha^{-1}(O) = O$ for all $\alpha$, in particular $\phi(j^1f)=j^1f$ for all $f \in C^\infty(N)$, since for any such $f$ there is $\alpha \in \cG$ such that $\alpha^{-1}(O) = j^1f$, and this forces $\phi = \id$ by lemma \ref{lemma_contactomorphism_fixing_jets_identity}. To see that there is $\alpha \in \cG$ such that $\alpha^{-1}(O) = j^1f$, let $h$ be a Hamiltonian on $T^*N$ given by cutting off $\pr^*(-f)$ outside a sufficiently large compact, where $\pr \fc T^*N \to N$ is the bundle projection. Then $\alpha = \phi_{\pi^*h}$ does the job, unless $\cG = \Cont_{0,c}(J^1N)$, in which case we need to cut off $\pi^*h$ outside a sufficiently large compact subset of $J^1N$. Here $\pi$ is the symplectic projection onto $T^*N$ whose domain ($J^1N$ or $T^*N \times S^1$) depends on the choice of the group $\cG$. The statements about the monotonicity of the identity maps are obvious. \qed
\end{prf}

Now we pass to
\begin{prf}[of proposition \ref{prop_EP_implies_spectral_partial_order}]
If $H,K$ are two contact Hamiltonians on $J^1N$ with $K \geq H$, then the isotopy $(\phi_H^t)^{-1}\phi_K^t$ is nonnegative, that is generated by a nonnegative contact Hamiltonian, and therefore $\ell_-(\phi^{-1}\psi) \geq 0$, as follows from lemma \ref{lemma_sp_nums_and_ct_Hamiltonian}. Since the notion of nonnegativity is conjugation-invariant, we also have $\ell_-(\alpha\phi^{-1}\psi\alpha^{-1}) \geq 0$ for any $\alpha \in \Cont_0(J^1N)$. This implies
$$\lceil\ell_+(\alpha\psi^{-1}\phi\alpha^{-1})\rceil = -\lfloor\ell_-(\alpha\phi^{-1}\psi\alpha^{-1})\rfloor \leq 0$$
which by definition means $\psi^{-1} \leq \phi^{-1}$, therefore $\phi \leq \psi$. \qed
\end{prf}

\subsection{Metrics}\label{section_metrics}

Here we prove theorem \ref{thm_metrics}.
\begin{prf}

Recall the definition of the metrics:
$$\rho_{\osc}(\phi,\psi) = \max\big\{\lceil \ell_+(\alpha\phi\psi^{-1}\alpha^{-1})\rceil - \lfloor\ell_-(\alpha\phi\psi^{-1}\alpha^{-1})\rfloor\,|\,\alpha\in\Cont_{0,c}(V\big)\}$$
and
$$\rho_{\sup}(\phi,\psi) = \max\big\{\max\big(|\lceil \ell_+(\alpha\phi\psi^{-1}\alpha^{-1})\rceil|,|\lfloor\ell_-(\alpha\phi\psi^{-1}\alpha^{-1})\rfloor|\big)\,|\,\alpha\in\Cont_{0,c}(V)\big\}\,.$$
Both these metrics in fact come from norms on $\Cont_{0,c}(V)$. Namely, let
$$\rho_{\osc}(\phi) = \max\big\{\lceil \ell_+(\alpha\phi\alpha^{-1})\rceil - \lfloor\ell_-(\alpha\phi\alpha^{-1})\rfloor\,|\,\alpha\in\Cont_{0,c}(V)\big\}$$
and
$$\rho_{\sup}(\phi) = \max\big\{\max\big(|\lceil \ell_+(\alpha\phi\alpha^{-1})\rceil|,|\lfloor\ell_-(\alpha\phi\alpha^{-1})\rfloor|\big)\,|\,\alpha\in\Cont_{0,c}(V)\big\}\,.$$
We will show in a moment that both these functions are finite, conjugation-invariant, nondegenerate, symmetric (that is they take the same value on $\phi$ and $\phi^{-1}$) and satisfy the triangle inequality; in other words they are norms on $\Cont_{0,c}(V)$. Since the metrics are defined in terms of the norms, it suffices to prove these properties for the norms and the claims for the metrics will follow.

Therefore let us prove the claimed properties of the norms $\rho_{\osc},\rho_{\sup} \fc \Cont_{0,c}(V) \to \R$. First of all, note that if $\phi = \phi_H$, where $H$ is a Hamiltonian, then point (iii) of proposition \ref{prop_properties_sp_nums_contact_periodic} implies 
$$\rho_{\osc}(\phi)\leq \left\lceil\int_0^1\max H_t\,dt\right\rceil - \left\lfloor\int_0^1\min H_t\,dt\right\rfloor$$
and
$$\rho_{\sup}(\phi)\leq \max \left(\left|\left\lceil\int_0^1\max H_t\,dt\right\rceil\right|, \left|\left\lfloor\int_0^1\min H_t\,dt\right\rfloor\right|\right)\,,$$
which proves that $\rho_{\osc},\rho_{\sup}$ are well-defined.

Conjugation invariance and symmetry follow from the definition. Let us prove nondegeneracy for $\rho_{\osc}$. Let $\phi$ be such that $\rho_{\osc}(\phi) = 0$, that is $\lceil\ell_+(\alpha\phi\alpha^{-1})\rceil = \lfloor\ell_-(\alpha\phi\alpha^{-1})\rfloor$ for all $\alpha$. It follows that $\ell_+(\alpha\phi\alpha^{-1}) = \ell_-(\alpha\phi\alpha^{-1})$, which by the definition of $\ell_\pm$ and point (iv) of proposition \ref{prop_properties_sp_nums_Leg} implies $\widetilde\alpha\widetilde\phi\widetilde\alpha^{-1}(O) = j^1a(\alpha)$, where $a(\alpha)$ is a number depending on $\alpha$. For $f \in C^\infty(N)$ let $T_f \fc J^1N \to J^1N$ be defined by $T_f(q,p,z)=(q,p-d_qf,z+f(q))$. This is the contactomorphism generated by the contact Hamiltonian $(q,p,z) \mapsto f(q)$. It follows (after cutting $T_f$ suitably outside a large compact) that
$$T_{-f}\widetilde\phi(T_f(O))=j^1a(f)\,,$$
therefore
$$ \widetilde\phi(j^1f)=j^1(f+a(f))\,.$$
Lemma \ref{lemma_contactomorphism_fixing_jets_identity} now shows that $\widetilde\phi=\id$, therefore $\phi=\id$. For $\rho_{\sup}$ it suffices to note that $\rho_{\sup}(\phi)=0$ implies $\rho_{\osc}(\phi)=0$.

It remains to prove that the above functions satisfy the triangle inequality. For $\rho_{\osc}$ this follows from the triangle inequalities for spectral numbers, point (ii) of proposition \ref{prop_properties_sp_nums_contact_periodic}. For $\rho_{\sup}$ the proof of this fact is a little tedious. We present here the spirit of the argument, leaving the details to the reader.

Consider the following special case: $|\lceil \ell_+(\phi)\rceil| \geq |\lfloor \ell_-(\phi)\rfloor|$ and same for $\psi$. We then have
\begin{align*}
-|\lceil\ell_+(\phi)\rceil| - |\lceil \ell_+(\psi)\rceil| &\leq \lfloor\ell_-(\phi)\rfloor + \lfloor \ell_-(\psi)\rfloor\\
&\leq \lfloor \ell_-(\phi\psi)\rfloor \\
&\leq \lceil\ell_+(\phi\psi)\rceil\\
&\leq \lceil\ell_+(\phi)\rceil+\lceil\ell_+(\psi)\rceil\\
&\leq |\lceil\ell_+(\phi)\rceil|+|\lceil\ell_+(\psi)\rceil|\,.
\end{align*}
This shows that under the above assumption we obtained
$$\max\big(|\lfloor \ell_-(\phi\psi)\rfloor|,|\lceil\ell_+(\phi\psi)\rceil|\big) \leq \max\big(|\lfloor \ell_-(\phi)\rfloor|,|\lceil\ell_+(\phi)\rceil|\big) +\max\big(|\lfloor \ell_-(\psi)\rfloor|,|\lceil\ell_+(\psi)\rceil|\big)\,.$$
Similar considerations show that this inequality holds for any $\phi,\psi$. We now conjugate by $\alpha$ throughout and take the maximum over $\alpha$ to obtain the desired triangle inequality.

Next we prove that $(\Cont_{0,c}(V),\preceq,\rho_{\sup})$ is an ordered metric space. It suffices to prove that if $\id \preceq \phi \preceq \psi$ then $\rho_{\sup}(\phi) \leq \rho_{\sup}(\psi)$. We have
$$\lceil\ell_+(\phi)\rceil = \lceil\ell_+(\phi\psi^{-1}\psi)\rceil\leq \underbrace{\lceil\ell_+(\phi\psi^{-1})\rceil}_{\phi\preceq\psi\Rightarrow\leq 0} + \lceil\ell_+(\psi)\rceil \leq \lceil\ell_+(\psi)\rceil\,.$$
Similarly, $\lfloor\ell_-(\phi)\rfloor \leq \lfloor\ell_-(\psi)\rfloor$. Next, we have
$$-\lfloor\ell_-(\phi)\rfloor=\lceil\ell_+(\phi^{-1})\rceil = \lceil\ell_+(\psi\psi^{-1}\phi^{-1})\rceil\leq \underbrace{\lceil\ell_+(\psi^{-1}\phi^{-1})\rceil}_{\psi^{-1}\preceq\phi\Rightarrow\leq 0} + \lceil\ell_+(\psi)\rceil \leq \lceil\ell_+(\psi)\rceil\,,$$
and analogously
$$-\lceil\ell_+(\phi)\rceil = \lfloor\ell_-(\phi^{-1})\rfloor \leq \lfloor\ell_-(\psi)\rfloor\,,$$
and all these inequalities imply together
$$-\lfloor\ell_-(\psi)\rfloor \leq \lceil\ell_+(\phi)\rceil \leq \lceil\ell_+(\psi)\rceil\text{ and } -\lceil\ell_+(\psi)\rceil\leq \lfloor\ell_-(\phi)\rfloor \leq \lfloor\ell_-(\psi)\rfloor\,,$$
therefore
$$\max\big(|\lceil\ell_+(\phi)\rceil|,|\lfloor\ell_-(\phi)\rfloor|\big) \leq \max\big(|\lceil\ell_+(\psi)\rceil|,|\lfloor\ell_-(\psi)\rfloor|\big)\,.$$
Now we can conjugate $\phi$ and $\psi$ by $\alpha$ throughout, and taking the maximum over $\alpha$, we finally obtain the desired.

Let us prove the last assertion of the theorem, namely that there is an order-preserving isometric embedding $(\Z,\leq,|\cdot - \cdot|) \hookrightarrow (\Cont_{0,c}(V), \preceq, \rho_{\sup})$. Indeed, let $h \fc T^*N \to [0,1]$ be an autonomous Hamiltonian taking the value $1$ on the zero section and let $H = \pi^*h$. Then the map $k \mapsto \phi_{kH}$ does the job.\qed
\end{prf}

\subsection{Contact rigidity}\label{section_cont_rigidity}

In this subsection we homogenize the spectral number $\ell_+$ and apply the result to prove the contact rigidity announced above. For $\phi \in \Cont_{0,c}(V)$ put
$$\nu(\phi) = \lim_{k \to \infty} \frac{\ell_+(\phi^k)}{k}\,.$$
The limit exists because the sequence $\lceil\ell_+(\phi^k)\rceil$ is subadditive. The function $\nu \fc \Cont_{0,c}(V) \to \R$ thus defined has some nice properties.
\begin{prop}
\begin{enumerate}
\item $\nu(\phi^k) = k\nu(\phi)$ for $k \geq 0$;
\item $\nu$ is conjugation-invariant;
\item $\nu(\phi\psi)\leq\nu(\phi)+\nu(\psi)$ if $\phi,\psi$ commute;
\item if $U \subset T^*N$ is displaceable by a Hamiltonian diffeomorphism then $\nu(\phi) = 0$ for any $\phi$ generated by a contact Hamiltonian with support in $U \times S^1$; more generally, if $\phi$ is such and $\psi$ commutes with $\phi$ then $\nu(\phi\psi) = \nu(\psi)$;
\item if $H$ is a contact Hamiltonian supported away from the zero section then $\nu(\phi_H) = 0$;
\item if $H$ is a contact Hamiltonian such that $H\geq c$ ($\leq c$, $=c$) on the zero wall, where $c\in\R$, then $\nu(\phi_H)\geq c$ ($\leq c$, $=c$);
\item if $h$ is a Hamiltonian on $T^*N$ and $H$ is its prequantization, then $\nu(\phi_H) = \mu(\phi_h)$, where $\mu$ is the partial quasi-morphism on the Hamiltonian group defined in \cite{Monzner_Vichery_Zapolsky_partial_qms_qss_cot_bundles}.
\end{enumerate}
\end{prop}
\begin{prf}
(i) follows directly from the definition.

(ii) We have
$$\lceil \ell_+(\alpha\phi\alpha^{-1}) \rceil \leq \lceil \ell_+(\alpha) \rceil + \lceil \ell_+(\alpha^{-1}) \rceil + \lceil \ell_+(\phi) \rceil\,,$$
and analogously
$$\lceil \ell_+(\phi) \rceil \leq \lceil \ell_+(\alpha) \rceil + \lceil \ell_+(\alpha^{-1}) \rceil + \lceil \ell_+(\alpha\phi\alpha^{-1}) \rceil\,,$$
which together imply
$$|\lceil \ell_+(\alpha\phi\alpha^{-1})\rceil - \lceil \ell_+(\phi)\rceil| \leq \lceil \ell_+(\alpha) \rceil + \lceil \ell_+(\alpha^{-1}) \rceil \,.$$
Taking now the $k$-th power of $\phi$ we obtain
$$|\lceil \ell_+(\alpha\phi^k\alpha^{-1})\rceil - \lceil \ell_+(\phi^k)\rceil| \leq \lceil \ell_+(\alpha) \rceil + \lceil \ell_+(\alpha^{-1}) \rceil \,,$$
and upon homogenization with respect to $k$ we obtain the assertion.

(iii) We have
$$\lceil \ell_+((\phi\psi)^k)\rceil = \lceil \ell_+(\phi^k\psi^k)\rceil \leq \lceil \ell_+(\phi^k)\rceil + \lceil \ell_+(\psi^k)\rceil\,.$$
Now homogenize.

(iv) Let $H$ be a contact Hamiltonian supported in $U \times S^1$. There are Hamiltonians $\underline{h}$ and $\ol h$ on $T^*N$ with supports in $U$ such that
$$\pi^*\underline h \leq H \leq \pi^*\ol h\,.$$
This implies, by the above properties of spectral numbers, that
$$\ell^{T^*N}_-(\phi_{\underline h}) \leq \ell_-(\phi_H) \leq \ell_+(\phi_H) \leq \ell^{T^*N}_+(\phi_{\ol h})\,.$$
It is proved in \cite{Monzner_Vichery_Zapolsky_partial_qms_qss_cot_bundles} that if a Hamiltonian $h$ on $T^*N$ has support in a displaceable subset $U$ then its spectral invariants are all contained in $[-e(U),e(U)]$ where $e(U)$ is a certain constant depending on $U$, the so-called spectral displacement energy of $U$. This implies that
$$-e(U) \leq \ell_-(\phi_H)\leq \ell_+(\phi_H) \leq e(U)\,.$$
Note that although the Hamiltonians $\underline h, \ol h$ depend on $H$, the resulting bound is in fact independent of them. Therefore upon homogenization we obtain, as required:
$$\nu(\phi_H) = 0\,.$$
The more general assertion follows from this and point (iii).

(v) If a contact Hamiltonian is supported away from the zero section, then so is its contact vector field, therefore the contact isotopy it generates preserves the zero section and thus the associated spectral invariants are all zero, and then so is the value of $\nu$ on the time-$1$ flow of such a Hamiltonian.

(vi) First note that if $H=c$ on the zero wall, then $\widetilde{\phi_H^t}(O)=j^1(ct)$, which is a Legendrian whose spectral invariants all equal $ct$. Upon homogenization we obtain $\nu(\phi_H) = c$. If, for instance $H \geq c$ on the zero wall, we can find another Hamiltonian $K$ which equals $c$ on the zero wall and such that $H \geq K$ everywhere. It follows from corollary \ref{coroll_comparison} that $\ell_+(\phi_H^k) \geq \ell_+(\phi_K^k)$ for all $k > 0$; now homogenize.

(vii) follows from point (vi) of proposition \ref{prop_properties_sp_nums_Leg}. Indeed, $\mu$ is the homogenization of the Lagrangian spectral invariant. \qed
\end{prf}

We can now prove theorems \ref{thm_rigidity_superheavy} and \ref{thm_rigidity_products}.

\begin{prf}[of theorem \ref{thm_rigidity_superheavy}]Assume by contradiction that there is $\alpha \in \Cont_{0,c}(V)$ such that $\alpha(X \times S^1) \cap O = \varnothing$. Then there is a neighborhood $W$ of $X$ such that $\alpha(W \times S^1) \cap O = \varnothing$. Let $h$ be a smooth function supported in $W$ and taking the value $1$ on $X$, and let $H$ be its prequantization. Then $\nu(\phi_H) = \mu(\phi_h)=\zeta(h) = 1$ since $X$ is assumed to be superheavy. On the other hand, $\alpha\phi_H^t\alpha^{-1}$ is an isotopy supported away from $O$. Therefore
$$0 = \nu(\alpha\phi_H\alpha^{-1})=\nu(\phi_H) = 1\,,$$
and this contradiction proves the claim. \qed
\end{prf}

\begin{prf}[of theorem \ref{thm_rigidity_products}]We need to prove that the time-$1$ map of a Hamiltonian $H$ bounded away from zero cannot be the product $\phi_1\dots\phi_k$. Indeed, we have $\nu(\prod_j\phi_j)=0$. On the other hand, if $H \geq c > 0$ along the zero wall, for instance, then $\ell_+(\phi_H^k) \geq kc$, which implies $\nu(\phi_H) \geq 1$. \qed
\end{prf}

%Recall that there is a canonical injection $\Ham(T^*N) \to \cG$. The order $\leq$ induces a partial order on $\Ham(T^*N)$. There is a natural normal cone in $\Ham(T^*N)$ consisting of all the time-$1$ maps of flows of nonnegative Hamiltonians. The above discussion shows that the relation induced by this normal cone is in fact a partial order; denote this order by $\leq_{\text{\rm EP}}$ as well; then the identity map $(\Ham(T^*N),\leq_{\text{\rm EP}}) \to (\Ham(T^*N),\leq)$ is monotone.

%\bibliographystyle{plain}
\bibliographystyle{abbrv}
\bibliography{1}

\end{document}